# Matrix Recovery using Split Bregman

Anupriya Gogna  
Student Member, IEEE

Ankita Shukla  
Student Member, IEEE

Angshul Majumdar  
Member, IEEE

ECE Department, Indraprastha Institute of Information Technology-Delhi, INDIA

*Abstract*— In this paper we address the problem of recovering a matrix, with inherent low rank structure, from its lower dimensional projections. This problem is frequently encountered in wide range of areas including pattern recognition, wireless sensor networks, control systems, recommender systems, image/video reconstruction etc. Both in theory and practice, the most optimal way to solve the low rank matrix recovery problem is via nuclear norm minimization. In this paper, we propose a Split Bregman algorithm for nuclear norm minimization. The use of Bregman technique improves the convergence speed of our algorithm and gives a higher success rate. Also, the accuracy of reconstruction is much better even for cases where small number of linear measurements are available. Our claim is supported by empirical results obtained using our algorithm and its comparison to other existing methods for matrix recovery. The algorithms are compared on the basis of NMSE, execution time and success rate for varying ranks and sampling ratios.

*Keywords—augmented lagrangian, bregman divergence, nuclear norm, rank minimization, split bregman*

## I. INTRODUCTION

Matrix recovery problem is commonly encountered in various areas of engineering and technology such as collaborative filtering [1], control system (for system identification) [2], wireless sensor networks (WSN) (for sensor localization) [3], computer vision (for inferring scene geometry) [4], data mining and pattern recognition etc. In the absence of any additional knowledge it is impossible to fill in the matrix given its subsampled version; as there can be infinitely many ways in which missing values can be filled. However, if the matrix to be recovered is (estimated to be approximately) low rank, the problem becomes feasible as one can search for the lowest rank matrix amongst all possible candidate solutions.

We give a few examples of low rank matrix recovery. In the field of collaborative filtering, we are required to find items of interest for a particular user, based on his/her ratings on other items and other user's ratings on similar items. This entails the task of filling an incompletely observed user-item rating matrix. The rating matrix can be considered to be low rank as the user preference for any item is assumed to be influenced by only a small number of factors. Similarly in the area of WSN, the task is to complete a partially filled distance matrix (wherein each element represents the distance between a pair of sensor nodes). It is known that such a matrix has a rank equal to the dimensionality of the space plus two. For such cases, matrix recovery problem becomes equivalent to obtaining the lowest rank solution consistent with the observed set of measurements.

So far we have discussed about matrix completion from partially observed entries. A somewhat different problem arises in dynamic MRI reconstruction. In MRI, the data is acquired in the Fourier frequency space. In dynamic MRI, the problem is to reconstruct the frames in time, given the sampled Frequency space. If we assume that the MRI frames to be reconstructed are stacked as columns of a matrix, the thus formed matrix is of low rank. This is because, the frames are correlated in time and hence the columns of the matrix are not independent from each other. Thus, dynamic MRI reconstruction turns out to be a low rank matrix recovery problem from sampled Fourier measurements.

In this work, we solve the problem of recovering a low rank matrix from its lower dimensional projections. Formally this is expressed as follows

$$\min_Z\ rank(Z)$$
$$\text{Subject to } Y = A(Z) \qquad (1)$$

Where, A is the linear projection operator (for matrix completion, A is a binary mask, for dynamic MRI, it is a Fourier mapping), Y is the acquired measurements and Z is the low rank matrix to be recovered.

Unfortunately, the aforesaid problem (1) is NP hard with doubly exponential complexity. It has been proven in [5, 6] that when certain conditions are met, the NP hard rank minimization (1) can be surrogated by its nearest convex envelope – the nuclear norm (2),

$$\min_Z\ \|Z\|_*$$
$$\text{Subject to } Y = A(Z) \qquad (2)$$

Where the nuclear norm $\|Z\|_*$ is the sum of singular values of the matrix.

The contribution of our work is the use of Split Bregman algorithm to solve (2) for obtaining the lowest rank solutions consistent with the given observation matrix Y. There are a handful of algorithms to solve the matrix completion problem, but there is hardly any algorithm to solve the general low rank matrix recovery problem from its lower dimensional projections. In this work, we develop a general algorithm for low rank matrix recovery using Split Bregman approach. It is shown that use of Split Bregman technique improves the accuracy of results as well as the speed of convergence. The rest of the paper is organized as follows. In section II, the previous work in the domain of matrix completion is discussed. Section III describes briefly Bregman iterations and Split Bregman algorithm. In section IV we introduce our proposed algorithm. Section V includes the experiments and results. Section VI contains the conclusions.

## II. REVIEW OF PREVIOUS STUDIES IN MATRIX COMPLETION

As mentioned before, there are a handful of algorithms to solve the matrix completion problem. It is not possible to review all of them; here we review two popular ones - Singular Value Thresholding (SVT) [7] and Fixed Point Continuation (FPC) [8]. Both solve the problem of the form (3)

$$\min_{Z} \ \|Z\|_* \quad \text{Subject to } A(Z) = b \quad (3)$$

Where, Z is the matrix to be recovered given the observed vector b and A is a linear operator. Matrix completion is a special case where A is a binary mask.

SVT uses linearized Bregman iterations for solving (3) using the following iterates

$$\begin{aligned} Z^k &= D_\tau(A^T(Y^{k-1})) \\ Y^k &= Y^{k-1} + \delta_k(b - A(Z^k)) \end{aligned} \quad (4)$$

Where, $Y^0 = 0 \in R^{m \times n}$ and $D_\tau$ is the (soft thresholding) shrinkage operator. $\delta_k$ is the step size at the $k^{th}$ iteration. However, SVT is well suited only for very low rank matrices and require a large number of iterations to converge. Thus the run times for SVT is very high for cases where only a few linear measurements are available.

FPC algorithm for matrix completion proposed in [8], uses operator splitting to recover a low rank matrix via nuclear norm minimization (5).

$$\begin{aligned} Y^k &= Z^k - \tau A^T(AZ^k - b) \\ Z^{k+1} &= S_\tau(Y^k) \end{aligned} \quad (5)$$

Where, $S_a(.)$ is the matrix shrinkage operator that performs soft thresholding on the singular values of the argument (matrix) by an amount a.

It should be noted that both SVT and FPC were developed to solve the matrix completion problem and not the matrix recovery problem. However, both of them can be modified to solve the matrix recovery problem as shown here.

## III. BACKGROUND FOR FORMULATION

This section presents a brief description of the background for our formulation. Here we discuss the concepts of Bregman divergence and Split Bregman method; the underlying fundamentals for our work.

Bregman iterative algorithms [9] are used to solve many constrained optimization problems, such as Basis Pursuit and TV denoising, because of its better convergence behavior and improved accuracy of reconstruction. Bregman distance forms the basis of formulation of these algorithms. For a convex function E: X→R, where u, v ∈ X and p belongs to the set of sub gradient of the function, Bregman distance is given by $D_E^p$.

$$D_E^p(u,v) = E(u) - E(v) - \langle p, u-v \rangle \quad (6)$$

Bregman distance is a non-negative quantity but is not a conventional distance metric as it is not symmetric and does not follow triangular inequality. It is just a measure of separation of two points. Consider the following constrained optimization problem.

$$\min_{u} \ E(u) \quad \text{Subject to } H(u) = 0 \quad (7)$$

Where, E and H are convex functions defined on $R^n$ with H being differentiable. The corresponding unconstrained formulation can be written as

$$\min_{u} \ E(u) + \lambda H(u) \quad (8)$$

Equation (8) is solved iteratively by Bregman iterative algorithm in the following steps

$$\begin{aligned} u^{k+1} &= \min_{u} D_E^p(u, u^k) + \lambda H(u) \\ p^{k+1} &= p^k - \nabla H(u^{k+1}) \end{aligned} \quad (9)$$

Bregman provides excellent convergence results because H monotonically decreases with every iteration and converges to the minimizer of the function. An advancement over Bregman iterations is the linearized Bregman algorithm which makes computation efficient by combining the minimizing of unconstrained problem and Bregman update into a single step; which can be solved exactly. However, linearized Bregman cannot solve problems involving multiple L1 regularization terms. To handle this issue Split Bregman method [10] was proposed. Consider the objective function given below where $\Phi(u)$ and $H(u)$ are convex and H is differentiable.

$$\min_{u} \ |\Phi(u)|_1 + H(u) \quad (10)$$

The underlying concept of Split Bregman is decomposition of L1 and L2 terms such that they form different sub problems which can be solved way easily than the original composite objective function. Rewriting equation (10) by letting $d = \Phi(u)$ the constrained formulation becomes

$$\min |d|_1 + H(u) \quad \text{Subject to } d = \Phi(u) \quad (11)$$

The unconstrained equivalent of (11) is obtained by adding a penalization function to the problem as in (12).

$$\min_{u} \ |d|_1 + H(u) + \frac{\lambda}{2} \|d - \Phi(u)\|_2^2 \quad (12)$$

Now comparing this to the form in (7), we consider $E(u,d) = |d|_1 + H(u)$,

$$(u^{k+1}, d^{k+1}) = \min_u D_E^p(u, u^k, d, d^k) + \frac{\lambda}{2} \| d - \Phi(u) \|_2^2$$

$$p_u^{k+1} = p_u^k + \lambda (\nabla \Phi)^T (\Phi u^{k+1} - d^{k+1}) \quad (13)$$

$$p_d^{k+1} = p_d^k + \lambda (d^{k+1} - \Phi u^{k+1})$$

The 1st update step can be solved using ADMM (alternating direction method of multipliers) [11] by alternatively keeping one variable fixed and optimizing over the other. Simplified form of (13) is given below.

$$u^{k+1} = \min_u H(u) + \frac{\lambda}{2} \| d - \Phi(u) - b^k \|_2^2$$

$$d^{k+1} = \min_d |d|_1 + \frac{\lambda}{2} \| d - \Phi(u) - b^k \|_2^2 \quad (14)$$

$$b^{k+1} = b^k + (\Phi u^{k+1} - d^{k+1})$$

Since H(u) is smooth and differentiable everywhere, updation for u can be solved analytically. Solution for d is nothing but the solution for synthesis prior formulation and is obtained directly by skrinkage (soft thresholding) operator. Last step is the updation of Bregman variable.

One important advantage that Split Bregman provides over other algorithms is that one can keep lambda (the regularization parameter) constant to a value that achieves fast convergence.

## IV. Algorithm Formulation

The general problem is to solve an underdetermined linear system of equations:

$$y_{N \times 1} = A(Z_{m \times n}) \quad (15)$$

Where $A: R^{m \times n} \to R^N, N < m \times n$

This has infinitely many solutions. We are interested in a low-rank solution. As mentioned before, one ideally needs to solve a rank minimization problem (1) in order to achieve this; but it has been theoretically proven that under certain conditions the same (low rank) solution can be achieved by nuclear norm minimization. In this work, we proposed to derive a new algorithm for nuclear norm minimization.

So far, we discussed the noiseless scenario. In all practical situations, the system will be corrupted by some noise; and we will have (16) instead,

$$Y = A(Z) + \eta, \quad \eta \sim N(0, \sigma^2) \quad (16)$$

To recover Z from (16) one ideally needs to solve (17) [12]:

$$\min_Z \| Z \|_* \text{ subject to } \| Y - A(Z) \|_{fro}^2 \leq \varepsilon \quad (17)$$

Where 'fro' denotes the Frobenius norm and $\varepsilon = N\sigma^2$.

Solving the constrained problem directly is difficult, so we propose to solve its unconstrained Lagrangian (18) instead. The constrained and the unconstrained counterparts are the same for correct choice of $\lambda_1$ given $\varepsilon$. In this work, we assume that $\lambda_1$ is known.

$$\min_Z \frac{1}{2} \| y - Az \|_2^2 + \lambda_1 \| Z \|_* \quad (18)$$

Where, A is a projection operator, z is the vectorized (column concatenated) form of matrix Z i.e. $A(Z) \equiv Az$, where for "Az" formulation, A is a block diagonal matrix such that each block of A acts on a column of matrix Z. In this section we derive an algorithm to solve (18). Split Bregman and ADMM techniques are adopted to ensure better accuracy and faster convergence. We reformulate (18) as (19) by introducing a proxy variable W.

$$\min_{Z,W} \frac{1}{2} \| y - Az \|_2^2 + \lambda_1 \| W \|_*$$
$$\text{Subject to } Z = W \quad (19)$$

The above formulation can be converted into unconstrained convex optimization problem (20) by use of augmented Lagrangian and Split Bregman techniques as discussed in section III.

$$\min_{Z,W} \left(\frac{1}{2}\right) \| y - Az \|_2^2 + \lambda_1 \| W \|_* + \left(\frac{\eta}{2}\right) \| W - Z - B \|_2^2 \quad (20)$$

Where, B is the Bregman (relaxation) variable.

The use of 2nd L2 term (augmented Lagrangian) improves the robustness of the algorithm as it eliminates the need for strictly reinforcing the equality constraint while simultaneously penalizing for any deviation. Use of Bregman variable makes sure that the value of $\lambda_1$ and $\eta$ can be chosen to optimize the convergence speed. Hence the speed of algorithm is dependent on how fast we can optimize each of the sub problems. Equation (20) can be split into two simpler sub problems which can be solved by alternatively fixing one variable and minimizing over the other (via ADMM).

Sub problem 1

$$\min_Z \left(\frac{1}{2}\right) \| y - Az \|_2^2 + \left(\frac{\eta}{2}\right) \| W - Z - B \|_2^2 \quad (21)$$

Sub problem 2

$$\min_W \lambda_1 \| W \|_* + \left(\frac{\eta}{2}\right) \| W - Z - B \|_2^2 \quad (22)$$

The first sub problem (21) is a simple least square problem (Tikhonov regularization) which can be solved easily using any gradient descent algorithm such as conjugate gradient. Sub problem 2 (22) can be solved iteratively by soft thresholding operation (23) on singular values of matrix W [13, 14].

$$Soft(t, u) = sign(t) \max(0, |t| - u) \quad (23)$$

$$W \leftarrow Soft(Z + B, 4\lambda_1 / \eta) \quad (24)$$

Between every consecutive iteration of the two sub problems, Bregman variable is updated as follows

$$B^{k+1} \leftarrow B^k - (W^{k+1} - Z^{k+1}) \quad (25)$$

Using Bregman variable causes initial few iterations to return over-regularized results. However updation of Bregman variable (25) as iterations proceed makes sure that any information that is not captured is added back. This scheme makes sure that the run times are lower and accuracy of recovery is improved. For solving the least square problem, we used lsqr [15]. Complete algorithm is given below.

---

Initialization: $B^0 = 1$; $W^0 = 0$; $\eta = 0.001$; $\lambda_1 = 0.001$

while $k < 500$ &&
  $abs(Obj\_func(k) - Obj\_func(k-1)) < 1E-7$

$Z^{k+1} = [A^T A + \eta I]^{-1}[\eta(W^k - B^k) + A^T y]$

$W^{k+1} \leftarrow \left( Z^{k+1} + B^k, \dfrac{4\lambda_1}{\eta} \right)$

// Update the Bregman Variable

$B^{k+1} \leftarrow B^k - \left( W^{k+1} - Z^{k+1} \right)$

end while

---

There are two stopping criterions. Iterations continue till the objective function converges; by convergence we mean that the difference between the objective functions between two successive iterations is very small (1E-7). The other stopping criterion is a limit on the maximum number of iterations. We have kept it to be 500.

## V. RESULTS

Experiments were carried out in two parts. In the first part, we experimented with synthetic datasets. In the second part, we experimented on actual collaborative filtering problem.

### A. Simulations on Synthetic data

For testing the algorithm, matrices of dimension $250 \times 250$ of varying rank (lower than the dimensionality of the matrix) were constructed. They were sub-sampled (at varying sampling ratios) using a binary (linear) operator A as in (2). For our algorithm, the values of regularization parameters $\lambda_1$ and $\eta$ were chosen to be 0.001. Bregman variable was initialized to a vector containing all ones. These values were chosen to yield the best convergence behavior. We compared our results against those obtained using FPC and SVT. For FPC the value of mu_final was taken to be 0.01 and tolerance was taken to be 1E-3. For SVT maximum number of iterations were taken to be 500.

Fig. 1-3 illustrate the success rates of FPC, SVT and MSB algorithms as a function of sampling ratio for different ranks. To compute the success rate, 100 independent runs of each algorithm were carried out. This was done for all combinations of sampling ratios and ranks. Success rate was computed as the number of attempts (out of 100) in which NMSE of less than 1E-3 was achieved. It can be seen that our proposed approach (MSB) always yields the best recovery; i.e. it can successfully recover matrices from low sampling ratios unlike the other two.

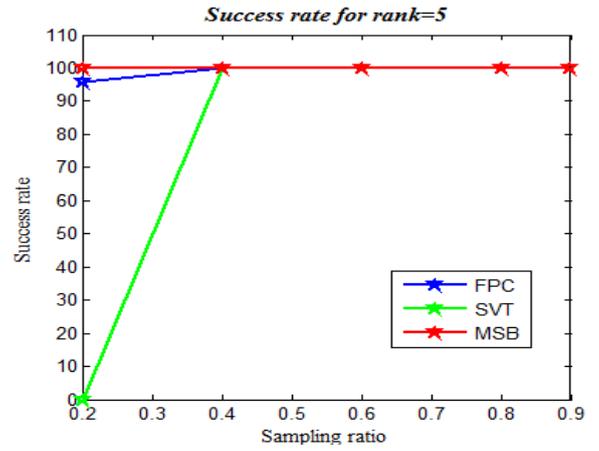

Fig. 1. Success rate for rank 5

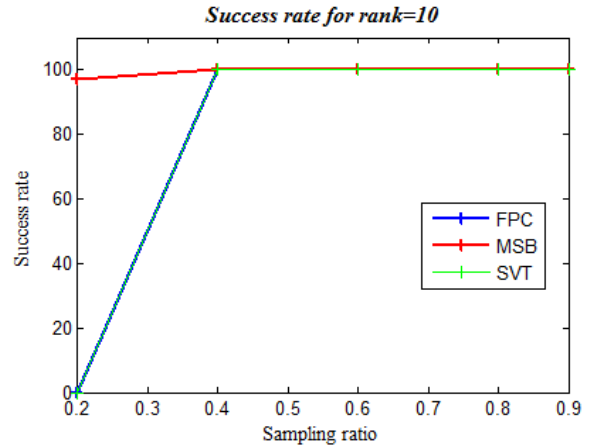

Fig. 2. Success rate for rank 10

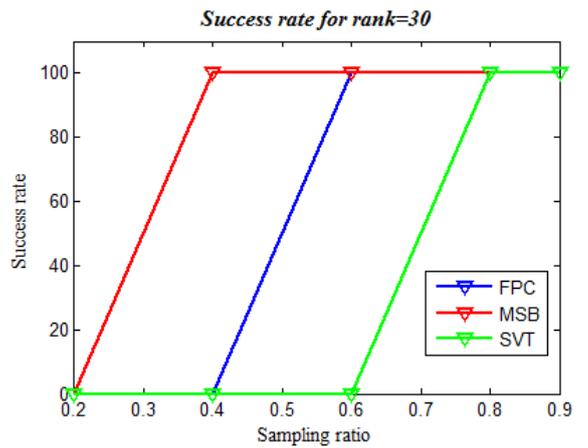

Fig. 3. Success rate for rank 30

The run times for all three algorithms are shown in table I. The table clearly shows that our algorithm converges faster to the optimum value. The execution time for MSB is much lower than that for both FPC and SVT especially for matrices with high rank and very few (available) sampled values. SVT has very poor run times for high rank matrices making it inefficient for use in such cases.

TABLE I. RUN TIMES FOR DIFFERENT ALGORITHMS

| Rank | Algorithm | Sampling Ratios | |
|---|---|---|---|
| | | *0.2 (sec)* | *0.4 (sec)* |
| 40 | FPC | 46.4 | 338.7 |
| | SVT | 447.2 | 372.6 |
| | MSB | 32.8 | 34.1 |
| 30 | FPC | 38.8 | 189.7 |
| | SVT | 406 | 179.3 |
| | MSB | 26.5 | 44.2 |
| 20 | FPC | 34.12 | 21.9 |
| | SVT | 342 | 31.7 |
| | MSB | 20.7 | 63.3 |

The success rate gives an overall picture. In order to assess the recovery accuracy, the success rate graphs are not enough. We need to observe the reconstruction errors. For this purpose, the average (Normalized Mean Squared Error) NMSE for all three algorithms (MSB, SVT and FPC), computed over 100 independent runs of each algorithm, is shown in table II. It is evident from the numerical values that our algorithm provides much lower NMSE, in all the cases, as compared to both SVT and FPC.

TABLE II. NMSE FOR DIFFERENT ALGORITHMS

| Rank | Algorithm | Sampling Ratios | | | |
|---|---|---|---|---|---|
| | | *0.2* | *0.4* | *0.6* | *0.8* |
| 5 | FPC | 5.1E-4 | 2.2E-4 | 1E-4 | 8.5E-5 |
| | SVT | 1.6E-3 | 2.1E-4 | 1.1E-4 | 1.01E-4 |
| | MSB | **6.1E-6** | **4.15E-7** | **1.06E-7** | **5.27E-8** |
| 10 | FPC | 1.67E-2 | 2E-4 | 8.99E-5 | 6.31E-5 |
| | SVT | 8.8E-3 | 6E-4 | 1E-4 | 1E-4 |
| | MSB | **2.45E-4** | **4.58E-7** | **5.89E-8** | **5.43E-8** |
| 20 | FPC | 4.19E-2 | 2.1E-4 | 7.31E-5 | 4.74E-5 |
| | SVT | 7.19E-2 | 2.3E-3 | 4.3E-4 | 1.3E-4 |
| | MSB | **1.81E-2** | **7.09E-6** | **6.96E-6** | **1.28E-6** |
| 30 | FPC | 4.39E-2 | 4.3E-3 | 6.98E-5 | 4.11E-5 |
| | SVT | 8.06E-2 | 1.18E-2 | 1.1E-3 | 2.1E-4 |
| | MSB | **4.33E-2** | **7.42E-4** | **9.43E-6** | **1.80E-6** |
| 40 | FPC | 4.4E-2 | 1.5E-2 | 7.44E-5 | 3.04E-5 |
| | SVT | 8.5E-2 | 2.78E-2 | 2.1E-3 | 1.6E-4 |
| | MSB | **4.1E-2** | **4.4E-3** | **3.58E-5** | **1.78E-5** |

*B. Collaborative Filtering*

Matrix completion, a special case of matrix recovery is required extensively in Collaborative Filtering (CF). In CF we have a sparsely filled rating matrix, which contains the ratings given by each user to some of the available items (e.g. Movies in the movielens data set). The task is to find the user's liking/disliking for other items (movies) based on the ratings explicitly provided by him/her for a few items. It is believed that each item can be characterized by certain features (for example genre, cast etc. for a movie) and users liking/disliking for an item is based on their affinity towards these features. Thus the underlying matrix structure is low rank as the number of features characterizing it are far less than the dimensionality of the rating matrix.

Here, we tested our algorithm on movielens (100K) data set [16]. This data set consists of 100,000 ratings (1-5) from 943 users on 1682 movies. Each user has rated at least 20 movies. We divided the data set into test and training data. Out of 100K values, 80K values were used as part of the training set and test set consisted of the remaining 20K values. Using training data set matrix recovery was performed and (Mean Absolute Error) MAE was computed between the test set data and the recovered values. MAE obtained using SVT and FPC is 1.41 and 0.82 respectively. Our proposed algorithm gives MAE of 0.78. All values were computed over 50 independent runs.

The fame of collaborative filtering owes to the famous 1 million dollar Netflix prize. The prize was won by the group who brought down the relative error by 8.43%. In this work, we have improved upon FPC by reducing the error by 5%. This is significant improvement by collaborative filtering standards.

VI. CONCLUSION

In this paper we presented an algorithm for recovering a low rank matrix given it's under sampled measurements. Ideally this requires solving a rank minimization problem; but as rank minimization is NP hard, it is surrogated by nuclear norm minimization. Nuclear norm minimization is a convex problem which can be solved via semi-definite programming. However SDP is slow and over the years, various fast alternatives have been proposed.

Almost all of these studies were formulated around the low-rank matrix completion problem. There is hardly any off-the-shelf algorithm for the general matrix recovery problem. In this work we use Bregman iterations and Split Bregman method to derive an algorithm for matrix recovery. Our algorithm shows faster convergence and much smaller run times than existing methods compared against. The accuracy of recovered results, quantified in terms of NMSE, and the empirical success rate of our algorithm is also far superior to that obtained from the methods compared against.

In case of real world collaborative filtering problem also, our algorithm is able to achieve a reduction of 5% in MAE as compared to existing methods (FPC and SVT).